\documentclass[article,11pt]{amsart}
\usepackage{amsaddr}
\usepackage{amssymb, amsthm, amsfonts, amsgen, stmaryrd}
\usepackage{slashed}
\usepackage[english]{babel}
\usepackage{fullpage}
\usepackage[centertags]{amsmath}
\usepackage{indentfirst}
\usepackage{fancyhdr}
\usepackage[dvips]{graphicx}
\usepackage{psfrag}
\usepackage{enumerate}
\numberwithin{equation}{section}
\usepackage[latin1]{inputenc}

\usepackage{hyperref}
\hypersetup{
pdftitle={},%
pdfauthor={},%
pdfsubject={},%
pdfkeywords={},%
colorlinks=true,%
linkcolor=blue,%
citecolor=red,%
linktocpage=true,%
hyperfootnotes=true,%
pageanchor=true
}
\usepackage{color}

\title{Derivation based differential calculi for noncommutative algebras deforming a class of  three dimensional  spaces } 
\date{16 May 2018}
\author{Giuseppe Marmo, Patrizia Vitale and Alessandro Zampini}
\address{Dipartimento di Fisica ``E. Pancini'', Universit\`a di Napoli Federico II, \\
Via Cintia - 80126 Napoli, Italy
\\
INFN - Sezione di Napoli, Via Cintia - 80126 Napoli, Italy } 
\email{marmo@na.infn.it}
\email{patrizia.vitale@na.infn.it} 
 \email{azampini@na.infn.it}

\newcommand{\nn}{\nonumber}

\newcommand{\dd}{{\rm d}}

\newcommand{\figureheight}{8cm}
\newcommand{\putfig}[2]{\begin{figure}[htp]
        \special{isoscale c:/itex/texfig/#1.wmf, \the\hsize \figureheight}
        \vspace{\figureheight}
        \caption{#2}\label{fig:#1}
        \end{figure}}
\newcommand{\pictureheight}{4cm}
\newcommand{\putpicture}[2]{\begin{figure}[htp]
        \special{isoscale c:/itex/texfig/#1.wmf, \the\hsize \pictureheight}
        \vspace{\pictureheight}
        \caption{#2}\label{fig:#1}
        \end{figure}}

\newcommand{\beqa}{\begin{eqnarray}}
\newcommand{\eeqa}{\end{eqnarray}}
\newcommand{\beq}{\begin{equation}}
\newcommand{\eeq}{\end{equation}}

\newcommand{\del}{\partial}

\newcommand{\R}{{\mathbb{R}}}

\newcommand{\N}{{\mathbb{N}}}
\newcommand{\C}{\mathbb{C}}

\begin{document}

\thispagestyle{empty}

\begin{abstract}
We equip a family of algebras whose noncommutativity is of Lie type with a derivation based differential calculus obtained, upon suitably using both inner and outer derivations,  as a reduction of a redundant calculus over the Moyal four dimensional space.

\end{abstract}


\maketitle
\tableofcontents

\section*{Introduction}
\label{sec:intro}

This paper addresses the problem of introducing a differential calculus on the family of  (Lie type)  non commutative algebras introduced in \cite{sele}. They are realised as a deformation   
of 
 the algebra of smooth functions on a class of three-dimensional Poisson algebras related to coadjoint orbits of Lie algebras.  More specifically, we are mainly interested here in deepening the analysis, initiated more than ten years ago in \cite{dica},  of  differential calculi which are obtained from suitable  Lie algebras of inner and outer derivations.  Such derivations satisfy the standard Leibniz rule, so our approach differs from the one developped in \cite{meli1,meli2}, where a quantum phase space on a Lie algebra type non commutative space is defined by deforming the coproduct on a suitable Hopf algebra.

Such geometric structures are of interest within the noncommutative formalism for  field and gauge  theory. In the framework advocated here, these are defined in terms of noncommutative algebras of fields whose dynamics requires a  Dirac operator, a Laplacian,   and  properly defined gauge connections, which can be  formulated in terms of  the noncommutative analogue of the Koszul gauge  connection \cite{wallet09,cmw11,degmw}. 

 Let us mention here that other approaches to noncommutative field theory have been developped. The  Connes spectral action \cite{connes} (also see \cite{Lizzireview} and references therein for an updated review)  relies  on a  definition of the differential calculus based on the notion of spectral triples, while a quantum group gauge theory on quantum spaces (see \cite{castellani, oxford, istvan}) is based on the notion of covariant calculi (see \cite{woro1,woro2}).     
We shall not discuss further these approaches in the present article, since the algebras we consider come as  subalgebras of the well known Moyal four dimensional one.  

 Another widely studied approach which goes under the name of twist-deformed  field theory concerns noncommutative algebras whose  product is obtained  by composing the ordinary commutative product  with  a so-called twist operator (see  \cite{twist} and references therein). In such a case the differential calculus is obtained as a  twist of  the standard differential calculus, see \cite{aschieri}.  However,  Lie algebra based non commutative  products cannot be formulated in term of a twist, therefore the above deformation procedure does not apply. We shall not discuss twist-deformed algebra in this paper and refer to the literature for more details.  
 
 The paper is organized as follows. In section \ref{review} we review the definition of derivation based differential calculus and the Jordan-Schwinger map which realizes three dimensional Lie algebras as Moyal subalgebras. In section \ref{sec:semi}  we address the problem of defining a differential calculus for semisimple subalgebras, while in section \ref{sec:not} we consider non semisimple Lie algebras. We conclude  in section \ref{concl} with a short summary and by pointing out possible new directions.

\section{The general setting}
\label{review}
In order to provide  the setting recalled in the (quite long) title of the paper, 
 we start by describing what a derivation based differential calculus on a unital algebra is. Since the 
noncommutative algebras we shall consider are realized  as  subalgebras of the four dimensional  Moyal algebra, we present such algebra together with its unitizations, and then the Lie algebra type noncommutativity which deforms the classical three dimensional space.

\subsection{Derivation based differential calculi}
\label{subsec:dbdc}
Given an orientable  $N$-dimensional differentiable manifold $M$, it is well known that the differential calculus on it  is the differential graded algebra $(\dd, \Omega(M)=\oplus_{k=0}^N\Omega_k(M))$, with $\Omega_{k}(M)$ the set of $k$ exterior forms and $\dd:\Omega_k(M)\to\Omega_{k+1}(M)$ the (graded)  exterior differential with $\dd^2=0$. Notice that $\mathcal F(M)=\Omega_0(M)$ gives the (commutative) algebra of smooth functions on $M$. The $\mathcal F(M)$-bimodule of 1-forms is dual to the set of vector fields 
$\mathfrak{X}(M)$. The set $\mathfrak X(M)$ coincides canonically with the space of all derivations of $\mathcal F(M)$, the commutator   $[X_1,X_2]f=X_1(X_2f)-X_2(X_1f)$ (with $X_{1,2}\in\mathfrak X(M)$ and $f\in\mathcal F(M)$) provides $\mathfrak X(M)$ with an infinite dimensional Lie algebra structure. 

When, following the approach of the Gelfand duality, the algebra $\mathcal F(M)$ is replaced by a suitable  non commutative    algebra $\mathcal A$, the problem of defining a differential calculus for it has been widely studied following different approaches, as we already mentioned in the Introduction. The approach we adopt in this paper is to explore how, given a (finite dimensional) Lie algebra of derivations acting $\mathcal A$, it is possible to dually define a $\mathcal A$-bimodule of forms and a whole  differential graded algebra that we interpret as a differential calculus on $\mathcal A$.

Assume indeed (see \cite{segal67, segal70, lama90, dv01})  that    $\mathfrak l$ is  a Lie algebra acting upon a unital associative algebra $\mathcal A$ by derivations, i.e. $\rho:\mathfrak l\,\to\,\mathrm{End}(\mathcal A)$ is a linear map with $[\rho(X_a),\rho(X_b)]=\rho([X_a,X_b])$ and 
 $\rho(X)(a_1a_2)\,=\,(\rho(X)a_1)a_2+a_1(\rho(X)a_2)$ for any $X,X_a,X_b\in\mathfrak l$ and $a_1,a_2$ in $\mathcal A$. Denote by 
$C_{\wedge}^n(\mathfrak l, \mathcal A)$  the set\footnote{Notice that the space of derivations for a given algebra $\mathcal A$ is a module only with respect to the centre $Z(\mathcal A)$ of the algebra.} of $Z(\mathcal A)$-multilinear alternating mappings $\omega\,:\,X_1\wedge\dots\wedge X_n\,\mapsto\,\omega(X_1, \dots, X_n)$ from $\mathfrak l^{\otimes n}$ to $\mathcal A$. On the graded vector space $C_{\wedge}(\mathfrak l, \mathcal A)=\oplus_{j=0}^{j=\mathrm {dim}\,\mathfrak l}C_{\wedge}^n(\mathfrak l, \mathcal A)$, with $C_{\wedge}^0(\mathfrak l, \mathcal A)=\mathcal A$, one can define a wedge product by (with $\omega\in C_{\wedge}^k(\mathfrak l, \mathcal A),\,\omega'\in C_{\wedge}^s(\mathfrak l, \mathcal A)$ and $X_j\,\in\,\mathfrak l$)
\beq
\label{wedp}
(\omega\wedge\omega')(X_1,\dots,X_{k+s})\,=\,\frac{1}{k!s!}\,\sum_{\sigma\in \mathcal S_{k+s}}(\mathrm{sign}(\sigma))\omega(X_{\sigma(1)},\dots,X_{\sigma(k)})\omega'(X_{\sigma(k+1)},\dots,X_{\sigma(k+s)})
\eeq
(where $\mathcal S_{k+s}$ is the set of permutations of $k+s$ elements) and 
the operator $\dd:C_{\wedge}^n(\mathfrak l, \mathcal A)\to C_{\wedge}^{n+1}(\mathfrak l, \mathcal A)$  by 
\begin{align}
(\dd\omega)(X_0, X_1,\dots, X_n)&=\,\sum_{k=0}^n(-1)^k\rho(X_k)(\omega(X_0,\dots,\hat{X}_k,\dots,X_n))\, \nn \\ &\qquad +\,\frac{1}{2}\,\sum_{r,s}(-1)^{k+s}\omega([X_r,X_s], X_0,\dots, \hat{X}_r,\dots,\hat{X}_s,\dots, X_n)
\label{ddef}
\end{align}
(with $\hat{X}_r$ denoting that the $r$-th term is omitted). Such operator  is easily proven to be a graded antiderivation with $\dd^2=0$, so $(C_{\wedge}(\mathfrak l, \mathcal A), \dd)$ is a graded differential algebra. Although the relations \eqref{wedp} and \eqref{ddef} are valid for both commutative and non commutative algebras, when the algebra $\mathcal A$ is not commutative one easily see that it is  in general $f_1\dd f_2\neq (\dd f_2)f_1$ and $\omega\wedge\omega'\neq-\omega'\wedge\omega$. One has indeed 
$(f_1\dd f_2)(X)\,=\,f_1\,(\rho(X)f_2)$ and $((\dd f_2)f_1)(X)\,=\,(\rho(X)f_2)f_1$, while $(\omega\wedge\omega')(X_1,X_2)\,=\,\omega(X_1)\omega'(X_2)-\omega(X_2)\omega'(X_1)$ and $(\omega'\wedge\omega)(X_1,X_2)\,=\,\omega'(X_1)\omega(X_2)-\omega'(X_2)\omega(X_1)$.
This exterior algebra is an example of a derivation based calculus, where the derivations come from the action of the Lie algebra $\mathfrak l$ upon the algebra $\mathcal A$. Such exterior algebra is often denoted by $\underline{\Omega}_{\mathfrak l}(\mathcal A)$, while its subset $\Omega_{\mathfrak l}(\mathcal A)\subset\underline{\Omega}_{\mathfrak l}(\mathcal A)$ is  defined as the smallest differential graded subalgebra of $\underline{\Omega}_{\mathfrak l}(\mathcal A)$ generated in degree $0$ by $\mathcal A$. By construction, every element in $\Omega^n_{\mathfrak l}(\mathcal A)$
can be written as a sum of $a_0\dd a_1\wedge\dots\wedge\dd a_n$ terms with $a_j\in\mathcal A$, while this is not necessary for elements in $\underline{\Omega}_{\mathfrak l}(\mathcal A)$.
This difference will be seen in some of the examples we shall describe in the following sections. It is easy to prove that, if the algebra $\mathcal A$ is the set of smooth functions on a paracompact manifold, then 
$\underline{\Omega}_{\mathfrak l}(\mathcal A)=\Omega_{\mathfrak l}(\mathcal A)$.

Upon the graded differential algebra $C_{\wedge}(\mathfrak l,\mathcal A)=\underline{\Omega}_{\mathfrak l}(\mathcal A)$ a contraction operator can be defined. If $X\in\mathfrak l$, then 
\beq
\label{contX}
(i_{X}\omega)(X_1,\dots,X_n)\,=\,\omega(X,X_1,\dots,X_n)
\eeq
(with $X_j\in\mathfrak l$) gives a degree $(-1)$ antiderivation from $C_{\wedge}^{n+1}(\mathfrak l, \mathcal A)\,\to\,C_{\wedge}^n(\mathfrak l,\mathcal A)$. The operator defined by $L_{X}=i_X\dd+\dd i_X$ is the degree zero Lie derivative along $X$, the set $(C_{\wedge}(\mathfrak l, \mathcal A)=\underline{\Omega}_{\mathfrak l}(\mathcal A), \dd, i_X, L_X=i_X\dd+\dd i_X)$ gives a Cartan calculus on 
$\mathcal A$ depending on  the Lie algebra $\mathfrak l$ of derivations.

\subsection{The Moyal algebra}
\label{subsec:moyal}
Consider  the symplectic  vector space $(\R^{2N}, \omega)$, where $\omega$ is the non degenerate closed (i.e. symplectic)  2-form on $\R^{2N}$ which is translationally invariant,  that we write as  $\omega=\dd q_a\wedge\dd p_a$ along a global Darboux coordinate system given by $(q_a, p_a)_{a=1,\ldots,N}$. The corresponding Moyal product (see \cite{marse, alg1, alg2})   reads (with $\theta>0$)  
\beq
\label{s1}
(f*g)(x)\,=\,\frac{1}{(\pi\theta)^{2N}}\int\int\dd u\dd v \,f(x+u)g(x+v)\,e^{-2i\omega^{-1}(u,v)/\theta}, \qquad\qquad \omega=\begin{pmatrix} 0 & 1_N \\ -1_N & 0 \end{pmatrix}
\eeq
 for $f,g\,\in\,{\mathcal S}(\R^{2N})$, i.e. the Schwartz space in $\R^{2N}$. The set 
\beq
\label{s2}
\mathcal A_{\theta}\,=\,(\mathcal S(\R^{2N}), *)
\eeq
is a non unital pre - $C^*$-algebra.  This algebra has a tracial property, with 
\beq
\label{s5}
\langle f\mid g\rangle\,=\,\frac{1}{(\pi\theta)^N}\,\int \dd x\,f(x)*g(x)\,=\,\frac{1}{(\pi\theta)^N}\,\int \dd x\,f(x)g(x),
\eeq
so it is possible to define a Moyal product on a larger set than $\mathcal S(\R^{2N})$. If $T\in\mathcal S'(\R^{2N})$ (with $\mathcal S'(\R^{2N})$ the space of continuous linear functional on $\mathcal S(\R^{2N})$, i.e. the space of tempered distributions), its action (evaluation)  upon $g\,\in\,\mathcal S(\R^{2N})$ is denoted  by $\langle T\mid g\rangle\,\in\,\C$. With $f\,\in\,\mathcal S(\R^{2N})$ one defines $T*f$ and $f*T$ in $\mathcal S'(\R^{2N})$ by 
\begin{align}
&
\langle T*f\mid g\rangle\,=\,\langle T\mid f*g\rangle, \nn \\ & \langle f*T\mid g\rangle\,=\,\langle T\mid g*f\rangle.
\label{s6}
\end{align}
One considers the space of left and right multipliers
\begin{align}
&\mathcal M_{L}^{\theta}\,=\,\{T\,\in\,\mathcal S'(\R^{2N})\,:\,T*f\,\in\,\mathcal S(\R^{2N})\,\forall \,f\,\in\,\mathcal S(\R^{2N})\} \nn \\
&\mathcal M_{R}^{\theta}\,=\,\{T\,\in\,\mathcal S'(\R^{2N})\,:\,f*T\,\in\,\mathcal S(\R^{2N})\,\forall \,f\,\in\,\mathcal S(\R^{2N})\}; 
\label{s3}
\end{align}
the set $\mathcal M^{\theta}=\mathcal M_{L}^{\theta}\cap\mathcal M_{R}^{\theta}$ is a unital $*$-algebra\footnote{In \cite{alg1} one proves that, if $T\in\mathcal S'$, $R\in\mathcal M_{L}$ and $S\in\mathcal M_{R}$ then  $\langle R*T\mid f\rangle=\langle T\mid f*R\rangle$ and $\langle T*S\mid f\rangle=\langle T\mid S*f\rangle$ for any $f\in\mathcal S$. These give a meaningful definition of the Moyal product within the multiplier algebra. }. 
Such a unitization  is not unique: it turns to be indeed  the maximal compactification of $\mathcal A_{\theta}$ defined by duality. It contains polynomials, plane waves, Dirac's $\delta$ and its derivatives. Its classical limit 
\beq
\lim_{\theta\to 0}\,\mathcal M^{\theta}\,=\,\mathcal O_{M}
\label{s4}
\eeq
is the set of smooth functions of polynomial growth on $\R^{2N}$ in all derivatives.  If  $\{~,~\}$ denotes the Poisson bracket structure on $\R^{2N}$ corresponding to the  symplectic form $\omega$,  the Moyal product has, on a suitable subset of $\mathcal M^{\theta}$, the asymptotic expansion in $\theta$ given by
\beq
\label{pros1}
f*g\,\sim\,fg\,+\,\frac{i\theta}{2}\{f,g\}\,+\,\sum_{k=2}^{\infty}(\frac{i\theta}{2})^k\frac{1}{k!}D_{k}(f,g)\qquad\mathrm{as}\,\,\theta\,\,\to\,\,0
\eeq
with $D_k$ the $k$-th order bidifferential operator which, in the  easiest case $(\R^2, \omega=\dd q\wedge\dd p)$,  is written as
\beq
\label{trao}
D_{k}(f,g)\,=\,\frac{\del^kf}{\del q^k}\,\frac{\del^kg}{\del p^k}\,-\,\left(\begin{array}{c} k \\ 1\end{array}\right)\,\frac{\del^kf}{\del^{k-1}q\del p}\,\frac{\del^k g}{\del^{k-1}p\del q}\,+\,\dots\,+\,(-1)^k\,\frac{\del^k f}{\del p^k}\,\frac{\del^kg}{\del q^k}. 
\eeq
If $f,g\in\mathcal M^{\theta}$,  their commutator is defined as 
\beq
\label{dcs}
[f,g]_{\theta}=f*g-g*f.
\eeq
Since  from \eqref{trao} we see  $D_{s}(f,g)=(-1)^sD_{s}(g,f)$,  we have that
\beq
\label{mobra}
[f,g]_{\theta}\,=\,i\theta\,\{f,g\}\,+\,\sum_{s=1}^{\infty}\frac{2}{(2s+1)!}\,\left(\frac{i\theta}{2}\right)^{2s+1}D_{2s+1}(f,g).
 \eeq
Upon denoting by 
$$
\mathcal P_{k}\,=\,\{q_1^{a_1}q_2^{a_2}p_1^{b_1}p_2^{b_2}\,:\,a_1+a_2+b_1+b_2=k\,\in\,\N\},
$$
i.e. the vector space of order $k$ polynomials in $\R^4$, the relation \eqref{mobra} allows to immediately see that  
 the elements in $\mathcal P_1$  fullfill the so called  canonical commutation relations
\begin{align}
&[q_a,q_b]_{\theta}=0, \nn \\ &[p_a,p_b]_{\theta}=0, \nn \\ & [q_a, p_b]_{\theta}=i\theta\delta_{ab},
\label{cr1}
\end{align}
while, if at least one among the elements  $f, \,g$ is in $S=\mathcal P_0\oplus\mathcal P_1\oplus\mathcal P_2$, it is
\beq
\label{cr}
[f,g]_{\theta}\,=\,i\theta\{f,g\}.
\eeq
This means  that $(S, \,\{~,~\})$ is a Poisson subalgebra of $\mathcal F(\R^4)$, while $(S, [~,~]_{\theta})$ is a Lie subalgebra in  $\mathcal M^{\theta}$ with respect to the $*$-product commutator \eqref{dcs}, indeed  isomorphic to a one dimensional central extension of  the Lie algebra $\mathfrak{isp}(4,\R)$  providing the infinitesimal generators of the   inhomogeneous  symplectic linear  group $\mathrm{ISP}(4,\R)$. Moreover, $S$ is the maximal Lie algebra acting upon both $\mathcal F(\R^4)$ and $\mathcal M^{\theta}$ in terms of derivations, via the operators 
\beq
\label{haop}
X_f\,:\,g\,\mapsto\,i\theta\{f,g\}\,=\,[f,g]_{\theta}.
\eeq


\subsection{Lie algebra type noncommutative spaces}
\label{subsec:lie}
We know from \cite{poipo} that any three dimensional Lie algebra $\mathfrak g$  with basis $(x_1,x_2, x_3)$ and commutator structure
\beq
\label{ls1}
[x_a,x_b]\,=\,f_{ab}^{\,\,\,\,c}x_c
\eeq
can be cast in the form
\begin{align}
&[x_1,x_2]=cx_3+hx_2, \nn \\
&[x_2,x_3]=ax_1, \nn \\
&[x_3,x_1]=bx_2-hx_3
\label{la3}
\end{align}
with real parameters $a,b,c,h$ such that $ah=0$, while from \cite{jsnap} we know that  a (classical) Jordan -- Schwinger map $\pi_{\mathfrak g}\,:\,\R^4\,\to\,\mathfrak g^*\sim\R^3$ can be defined such that 
\beq
\{\pi_{\mathfrak g}^*(x_a), \pi_{\mathfrak g}^*(x_b)\}\,=\,f_{ab}^{\,\,\,\,c}\pi_{\mathfrak g}^*(x_c).
\label{JSc}
\eeq
Upon noticing  that a classical Jordan -- Schwinger map with $\pi^*(x_a)$ ranging within $\mathcal P_1\oplus\mathcal P_2$ can be defined, a (noncommutative) version of it is introduced in \cite{sele} as a vector space inclusion $s_{\mathfrak g}\,:\,\mathfrak g^*\hookrightarrow\,\mathcal P_1\oplus \mathcal P_2$ such that
\beq
[s_{\mathfrak g}(x_a),s_{\mathfrak g}(x_b)]_{\theta}\,=\,i\theta f_{ab}^{\,\,\,\,c}\,s_{\mathfrak g}(x_c).
\label{qJS}
\eeq
On fixing  a three dimensional Lie algebra $\mathfrak g$, the Moyal product in $\R^4$ of functions of the variables $s_{\mathfrak g}(x_a)$ is proven to depend only on  the $s_{\mathfrak g}(x_a)$ variables, so there  exists a unital complex $*$-algebra $A_{\mathfrak g}\subset\mathcal M^{\theta}$ which is given  as the quotient
\beq
\label{Aqu}
A_{\mathfrak g}\,=\,[u_1,u_2,u_3] /I_{\mathfrak g},
\eeq
where (upon denoting $u_a=s_{\mathfrak g}(x_a)$) $I_{\mathfrak g}$ is the two-sided ideal generated by $\left[u_au_b-u_bu_a-i\theta f_{ab}^{\,\,\,\,c}u_c\right]$ (see also \cite{gutt83}: we are realizing the universal enveloping algebra $A_{\mathfrak g}$ as a subalgebra of $\mathcal M^{\theta}$). We list the maps $s_{\mathfrak g}$ from \cite{sele}, starting by those having a (quadratic, i.e. in $\mathcal P_2$)  Casimir function $C_{\mathfrak g}$. In this list we do not consider the Abelian Lie algebra.

\bigskip

\begin{enumerate}[(1)]
\item For $\mathfrak g=\mathfrak{su}(2)$, with 
$$
[x_a,x_b]=\varepsilon_{ab}^{\,\,\,\,c}x_c,
$$
 we have 
\beq
\label{jsu2}
u_1=\,\frac{1}{2}\,(q_1q_2+p_1p_2), \qquad u_2=\,\frac{1}{2}\,(q_1p_2-q_2p_1), \qquad u_3=\,\frac{1}{4}\,(q_1^2+p_1^2-q_2^2-p_2^2).
\eeq
The Casimir function is given by $C_{\mathfrak{su}(2)}=u_1^2+u_2^2+u_3^2$.  The identity 
\beq
 u_4^2=u_1^2+u_2^2+u_3^2, \qquad\qquad\mathrm{with}\quad u_4=\,\frac{1}{4}\,(q_1^2+p_1^2+q_2^2+p_2^2),
 \label{cacsup}
 \eeq
shows that $u_4$ is the quadratic Casimir for $A_{\mathfrak{su}(2)}$.

\bigskip

\item For $\mathfrak g=\mathfrak{sl}(2, \R)$ with 
\begin{align*}
&[x_1,x_2]=-x_3,  \\ &[x_2,x_3]=x_1,  \\ &[x_3,x_1]=-x_2, 
\end{align*}
 we have 
\beq
\label{jsl2}
u_1=\,\frac{1}{4}\,(q_1^2+p_1^2+q_2^2+p_2^2), \qquad u_2=\frac{1}{4}\,(q_1^2+q_2^2-p_1^2-p_2^2), \qquad u_3=\,\frac{1}{2}\,(q_1p_1+q_2p_2).
\eeq
The Casimir function is $C_{\mathfrak{sl}(2,\R)}=u_1^2-u_2^2-u_3^2$. Analogously to the previous case, the identity 
\beq
 u_4^2=u_1^2-u_2^2-u_3^2, \qquad\qquad\mathrm{with}\quad u_4=\,\frac{1}{2}\,(q_1p_2-q_2p_1), 
 \label{cacslp}
 \eeq
shows that $u_4$ is the quadratic Casimir for $A_{\mathfrak{sl}(2,\R)}$.

\bigskip

\item For $\mathfrak g=\mathfrak{e}(2)$, with 
\begin{align*}
& [x_1,x_2]=x_3, \\ &[x_2,x_3]=0, \\ &[x_3,x_1]=x_2,
\end{align*} 
we have
\beq
\label{je2}
u_1=\,q_1p_2-q_2p_1, \qquad u_2=q_1, \qquad u_3=q_2.
\eeq
The quadratic  Casimir function is $C_{\mathfrak{e}(2)}=(u_2^2+u_3^2)/2$.

\bigskip

\item For $\mathfrak g=\mathfrak{iso}(1,1)$ with 
\begin{align*}
&[x_1,x_2]=x_3, \\ &[x_2,x_3]=0, \\ &[x_3,x_1]=-x_2, \end{align*} 
 we have
\beq
\label{jiso11}
u_1=\,\frac{1}{2}\,(p_1^2+p_2^2-q_1^2-q_2^2), \qquad u_2=q_2+p_1, \qquad u_3=-q_1-p_2.
\eeq
The quadratic Casimir function is $C_{\mathfrak{iso}(1,1)}=(u_3^2-u_2^2)/2$.

\bigskip
\item for $\mathfrak g=\mathfrak h(1)$ (the Heseinberg-Weyl Lie algebra), with 
\begin{align*}
&[x_1,x_2]=x_3, \\ &[x_2,x_3]=0, \\ &[x_3,x_1]=0,
\end{align*}
 we have
\beq
\label{jh1}
u_1=q_1, \qquad u_2=q_2p_1, \qquad u_3=q_2.
\eeq
The quadratic Casimir function is $C_{\mathfrak{h}(1)}=u_3^2$.
\end{enumerate}
The three dimensional Lie algebras having no Casimir correspond to the case $a=0$ in \eqref{la3}. A  Jordan-Schwinger map for them is given by 
\beq
\label{sca}
u_1=-h(q_1p_1+q_2p_2)-cq_2p_1+bq_1p_2, \qquad u_2=q_1, \qquad u_3=q_2.
\eeq
Among the algebras we listed, the one corresponding to $\mathfrak g=\mathfrak{su}(2)$ (which is closely related to the  fuzzy sphere \cite{fs}, also see \cite{matrixLV} for a review on fuzzy spaces) has been intensively studied (see \cite{viwa13, gvw14, rovi12, gjw15, vitale14, kuvi15})
\footnote{Notice that such algebra had been considered  in \cite{shei}.} as an example for a space with a Lie algebra type noncommutativity, so to provide,  corresponding to a suitable differential calculus, a gauge action giving a renormalizable field theory. 

\bigskip

Further elaborating on our preliminary results in  \cite{dica}, in this paper we equip  the algebras $A_{\mathfrak g}\subset\mathcal M^{\theta}$ which   correspond to a Lie algebra $\mathfrak g$ having a quadratic Casimir function, with an exterior algebra and a differential calculus. These are realized as a suitable reduction  of the  calculus on $\mathcal M^{\theta}$ given by 
$(C_{\wedge}(\mathfrak{isp}(4,\R), \mathcal M^{\theta}), \dd)$, where the derivations are defined by \eqref{haop} (see \cite{cmw11, masson08, wallet09}). We leave the remaining cases to a further work.

It is well known that the Moyal algebra $\mathcal M^{\theta}$ is a normal space of distributions, and all its derivations are inner. For any algebra  
$ A_{\mathfrak g}\subset\mathcal M^{\theta}$, the union of  its inner and outer derivations close a Lie algebra $\tilde{\mathfrak g}$ with $\mathfrak g\subseteq\tilde{\mathfrak g}\subset\mathfrak{isp}(4,\R)$. The Lie algebra $\tilde{\mathfrak g}$  is seen to act  via inner derivations upon $\mathcal M^{\theta}$, and such action can be projected onto $ A_{\mathfrak g}$. The set $C_{\wedge}(\tilde{\mathfrak g}, A_{\mathfrak g})$ can be then described as a graded subalgebra of $C_{\wedge}(\mathfrak{isp}(4,\R), \mathcal M^{\theta})$, the corresponding calculus $(C_{\wedge}(\tilde{\mathfrak g}, A_{\mathfrak g}),\dd)$
as a reduction of $(C_{\wedge}(\mathfrak{isp}(4,\R), \mathcal M^{\theta}), \dd)$. 
 Moreover, the differential calculus that we define on  $ A_{\mathfrak g}$  turns out to have a frame, i.e. the exterior algebra is a free $ A_{\mathfrak g}$-bimodule: this gives  a way to study its cohomology.
 
Since the structure of the space of derivations for $A_{\mathfrak g}$ strongly depends on the Lie algebra $\mathfrak g$ being semisimple or not, the paper has two sections, which cover the two cases.  The subsection \ref{subsec:YM} moreover shows that, for the considered  differential calculi on $A_{\mathfrak g}$ with semisimple $\mathfrak g$, with respect to a natural  Cartan-Killing symmetric tensor on it, one can explicitly write a vector potential solving the Yang-Mills equations. This result mimics a result proven in \cite{HTS80, NT78} for suitable  classical fibrations.

\section{Differential calculus on $A_{\mathfrak g}$ for   semisimple Lie algebras 
}
\label{sec:semi}

Among the three dimensional Lie algebras considered in the introduction, the only semisimple ones are $\mathfrak{su}(2)$ and $\mathfrak{sl}(2,\R)$, which are the two real forms of $\mathfrak{sl}(2,\C)$. It follows that the corresponding algebras $A_{\mathfrak g}$ have no outer derivations, and that the Cartan-Killing form for them  is non degenerate.   Since the  semisimplicity of $\mathfrak g$ strongly affects the structure of the derivation based calculi on $A_{\mathfrak g}$,  we start by analysing such a case. 

\subsection{The case $\mathfrak g=\mathfrak{su}(2)$}\label{subsec:su2}  We begin upon analysing the classical (i.e. commutative) setting for this algebra, so we consider the symplectic vector space $(\R^4, \omega=\dd q_1\wedge\dd p_1+\dd q_2\wedge\dd p_2)$ equipped with the Euclidean metric $g=\dd q_a\otimes\dd q_a+\dd p_a\otimes\dd p_a$.  Given the elements $u_1,u_2,u_3$ as in \eqref{jsu2} with $\{u_a,u_b\}=\varepsilon_{abc}u_c$ and the element \eqref{cacsup}
\beq
\label{cacsu}
u_4=\,\frac{1}{4}\,(q_1^2+p_1^2+q_2^2+p_2^2),
\eeq
 the corresponding Hamiltonian vector fields $(X_1,X_2,X_3)$ coincide with the right invariant vector fields which are tangent to the $S^3$ orbits of the  standard action of ${\rm SU}(2)$ upon $\R^4\backslash\{0\}$, while $X_4$ is the infinitesimal generator for the ${\rm U}(1)$ action giving the well known Kustaanheimo-Stiefel fibration $\pi_{KS}:\R^4\backslash\{0\}\,\stackrel{{\rm U}(1)}{\longrightarrow}\,\R^3\backslash\{0\}$. The set
$$
A\,=\,\{f\,\in\,\mathcal F(\R^4\backslash\{0\})\,:\,X_4(f)=\{u_4,f\}=0\}
$$
is the basis algebra for the $\pi_{KS}$ fibration: one can equivalently say \cite{damv} that the derivations $(X_1,\dots, X_4)$ give  the projectable vector fields for such a fibration. As in \cite{sele}, we define an algebra $ \tilde A$ as
\beq
\label{Asu2}
\tilde{A}\,=\,\{f\,\in\,\mathcal M^{\theta}\,:\,[u_4,f]_{\theta}=0\}. 
\eeq
Such an algebra can be realised as a suitable completion of the polynomial algebra $P_{\mathfrak{su}(2)}\,=\,[u_1,u_2,u_3,u_4] /\mathfrak I_{\mathfrak{su}(2)}$ with $\mathfrak I_{\mathfrak{su}(2)}$ the two-sided ideal generated by $\left[u_au_b-u_bu_a-i\theta \varepsilon_{abc}u_c\right]$ and by $[u_4^2-(u_1^2+u_2^2+u_3^2)]$. Notice that $\tilde A$ extends the algebra $A_{\mathfrak{su}(2)}$ defined in \eqref{Aqu} for the $\mathfrak{su}(2)$ case since it contains the odd powers of $u_4$. 
The lines above show indeed that the algebra $\tilde A$ is a noncommutative deformation of the commutative algebra  $\mathcal F(\R^3\backslash\{0\})$. 

Within the classical setting, the rank of the space of derivations for $A$ is 3, since we can write 
\beq
\label{derad}
X_4\,=\,u_4^{-1}\sum_{j=1}^3u_jX_j
\eeq
on $\mathcal F(\R^4\backslash\{0\})$, while the derivations on $\mathcal M^{\theta}$ given by 
\beq
\label{desu2}
D_{\mu}(f)\,=\,[u_{\mu},f]_{\theta}, \qquad \mu=1,\dots,4
\eeq
are independent.  The elements $(u_1,\dots,u_4)$ from $\tilde A$ give a (one dimensional) central extension $\tilde{\mathfrak g}$ of the Lie algebra $\mathfrak g=\mathfrak{su}(2)$ with respect to the $*$-product commutator in $\mathcal M^{\theta}$, with clearly $\tilde{\mathfrak g}\subset\mathfrak{isp}(4,\R)$. The inclusion  $C_{\wedge}^1(\tilde{\mathfrak g}, \mathcal M^{\theta})\subset C_{\wedge}^1(\mathfrak{isp}(4,\R), \mathcal M^{\theta})$  has been studied in  
\cite{dica}. Using the definitions \eqref{ddef}, the elements  
\begin{align}
\alpha_1&=\,p_2*\dd q_1+p_1*\dd q_2-q_2*\dd p_1-q_1*\dd p_2, \nn \\
\alpha_2&=\,-q_2*\dd q_1+q_1*\dd q_2-p_2*\dd p_1+p_1*\dd p_2, \nn \\
\alpha_3&=\,p_1*\dd q_1-p_2*\dd q_2-q_1*\dd p_1+q_2*\dd p_2, \nn \\
\beta&=\,q_1*\dd q_1+q_2*\dd q_2+p_1*\dd p_1+p_2*\dd p_2
\label{vfor}
\end{align}
are seen to satisfy the identities
\begin{align}
\alpha_{j}(D_k)\,=\,-2i\theta\,\delta_{jk}u_4, &\qquad\qquad\alpha_j(D_4)\,=\,-2i\theta\,u_j \nn \\
\beta(D_k)\,=\,0,&\qquad\qquad\beta(D_4)\,=\,\theta^2
\label{vfor1}
\end{align}
 with $j,k=1,\dots,3$.
Upon defining 
\begin{align}
\omega_j&=\,\frac{i}{2\theta}\,\alpha_j\,-\,\frac{1}{\theta^2}\,u_k\beta, \nn \\
\omega_4&=\,\frac{1}{\theta^2}\,u_4\beta,
\label{vfor2}
\end{align}
we have 
$$\omega_{\mu}(D_{\sigma})=u_4\delta_{\mu\sigma},
$$
 with $\mu,\sigma=1,\dots,4$. Since $Z(\mathcal M^{\theta})=\C$, the elements $\{\omega_{\mu}\}_{\mu=1,\dots,4}$ are $\C$-linear maps from the Lie algebra $\tilde{\mathfrak g}$ of derivations to $\tilde A$.  Since  $u_4\in Z(\tilde A)$, we extend $\tilde{A}$ upon a localization, i.e. we define the element $u_4^{-1}$  via the relations $u_4^{-1}u_4=u_4u_4^{-1}=1$ and $u_4^{-1}u_k=u_ku_4^{-1}$ for $k=1,\dots,3$. Denoting the extended algebra by the same symbol, 
we consider the  vector  space $\mathcal D\simeq\tilde{\mathfrak g}$ as the tangent space to the noncommutative space described by the algebra $\tilde A$, and clearly have that 
the elements 
$$
\varphi_{\mu}\,=\,u_4^{-1}\omega_{\mu}
$$
  provide a basis for $\mathcal D^*$ which is $\C$-dual to $\{D_{\mu}\}_{\mu=1,\dots,4}$. We consider then the free $\tilde A$-bimodule $C^1_{\wedge}(\tilde{\mathfrak g}, \tilde{A})$ generated by the elements $\{\varphi_{\mu}\}$. For  $f\in \tilde A$ the action of the exterior derivative upon $\tilde A$ is given by 
\beq
\label{desu2}
\dd f\,=\,(D_{\mu}f)\varphi_{\mu}
\eeq
where $D_{\mu}f=[u_{\mu},f]_{\theta}$. The relations \eqref{wedp} and \eqref{ddef} can be used to define wedge products (so to have a graded algebra  $C_{\wedge}({\tilde{\mathfrak g}}, \tilde{A})=\oplus_{j=0}^{4}C^{j}_{\wedge}( {\tilde{\mathfrak g}}, \tilde A)$) and to extend the action of the exterior derivative. The same relations allow to prove also that, although is $\tilde{A}$ not commutative, we have 
\begin{align}
&f*\varphi_{\mu}\,=\,\varphi_{\mu}*f, \nn \\
&\varphi_{\mu}\wedge\varphi_{\sigma}=-\varphi_{\sigma}\wedge\varphi_{\mu}.
\label{eafr}
\end{align}
Using  such identities, it is also immediate to see that, for such  a four dimensional differential calculus on $\tilde A$, it is 
\begin{align}
&\dd \varphi_j\,=\,-\frac{1}{2}\,\varepsilon_{jkl}\varphi_k\wedge\varphi_l \qquad\qquad(j,k,l\in 1,\dots,3) \nn \\
&\dd \varphi_4\,=\,0.
\label{mcsu2}
\end{align}
This means that the Maurer-Cartan equation for the differential calculus depends on the structure constants of the Lie algebra $\tilde{\mathfrak g}$, and its cohomology is related to the Eilenberg-Chevalley cohomology for $\tilde{\mathfrak g}$. 

Notice that the elements $\varphi_a$ cannot be realised as $\sum_{a=1}^3f_a\dd u_a$, so 
$C_{\wedge}({\tilde{\mathfrak g}}, \tilde{A})$ 
 extends the differential calculus $(\Omega_{\tilde{\mathfrak g}}, \dd)$ given as the smallest  graded differential subalgebra of $C_{\wedge}(\tilde{\mathfrak g},\tilde{A })$ generated in degree $0$ by $\tilde A$ as described in the introduction.

\subsection{The case $\mathfrak g=\mathfrak{sl}(2,\R)$}
\label{subsec:sl2}
The noncommutative algebra $\tilde A$ can be defined, as in \cite{sele}, upon setting
\beq
\label{Asl2}
\tilde{A}\,=\,\{f\,\in\,\mathcal M^{\theta}\,:\,[u_4,f]_{\theta}=0\}, 
\eeq
where the elements $\{u_1,u_2,u_3,u_4\}\in\mathcal P_2\subset\mathcal M^{\theta}$ are given by \eqref{jsl2} and \eqref{cacslp}.

As in the previous case for $\mathfrak{su}(2)$, this algebra is a suitable completion of the polynomial algebra $P_{\mathfrak{sl}(2,\R)}\,=\,[u_1,u_2,u_3,u_4] /\mathfrak I_{\mathfrak{sl}(2,\R)}$ with $\mathfrak I_{\mathfrak{sl}(2,\R)}$ the two-sided ideal generated by $\left[u_au_b-u_bu_a-i\theta f_{ab}^{\,\,\,\,c}u_c\right]$ and by $u_4^2-(u_1^2-u_2^2-u_3^2)$ with respect to the Lie algebra structure constants of $\mathfrak{sl}(2,\R)$. Even in this case we see that $\tilde A$ extends the algebra $A_{\mathfrak{sl}(2,\R)}$ defined in \eqref{Aqu} for the $\mathfrak{sl}(2,\R)$ since it contains the odd powers of $u_4$. The algebra $\tilde A$ gives a noncommutative deformation of the commutative algebra $\mathcal F(\R^3\backslash\{\|u\|=0\})$ where we have denoted the three dimensional hyperbolic norm as $\|u\|^2=u_1^2-u_2^2-u_3^2$.  

The analysis for a derivation based  differential calculus for $\tilde A$ proceeds as in the $\mathfrak{su}(2)$ case, so we limit ourselves to report that the elements $(u_1,\dots,u_4)$ give a one dimensional central extension $\tilde{\mathfrak g}$ of $\mathfrak g=\mathfrak{sl}(2,\R)$ with respect to the $*$-commutator in $\mathcal M^{\theta}$. Given the basis of the tangent space $\mathcal D\simeq\tilde{\mathfrak{g}}$ with basis elements given by the derivations
\beq
\label{desl2}
D_{\mu}(f)\,=\,[u_{\mu},f]_{\theta},
\eeq
it is immediate to prove  that the elements 
\begin{align}
\alpha_1&=\,-p_2*\dd p_1+p_1*\dd p_2+q_2*\dd q_1-q_1*\dd q_2, \nn \\
\alpha_2&=\,q_2*\dd p_1-q_1*\dd p_2+p_2*\dd q_1-p_1*\dd q_2, \nn \\
\alpha_3&=\,-p_2*\dd p_1+p_1*\dd p_2-q_2*\dd q_1+q_1*\dd q_2, \nn \\
\beta&=\,p_2*\dd q_1-p_1*\dd q_2-q_2*\dd p_1+q_1*\dd p_2
\label{vforsl}
\end{align}
in $C^1_{\wedge}(\mathfrak{isp}(4,\R),\mathcal M^{\theta})$ allow to define the forms 
\begin{align}
\omega_j&=\,\frac{i}{2\theta}\,\alpha_j\,-\,\frac{1}{\theta^2}\,u_k\beta, \nn \\
\omega_4&=\,\frac{1}{\theta^2}\,u_4\beta,
\label{vforsl2}
\end{align}
that  give
$$\omega_{\mu}(D_{\sigma})=u_4\delta_{\mu\sigma},
$$
 with $a,b=1,\dots,4$. After the natural localization provided by adding  $\tilde A$ the element $u_4^{-1}$ we see that, as in the previous case, the elements $\varphi_a=u_4^{-1}\omega_a$ give a basis for $\mathcal D^*$ which is dual to the basis \eqref{desl2} for $\mathcal D$. Such a basis generates, as already described, the four dimensional calculus $(C_{\wedge}({\tilde{\mathfrak g}}, \tilde{A}),\dd)$  whose Maurer-Cartan relations are clearly given by 
\begin{align}
\dd\varphi_1&=\,-\,\dd\varphi_2\wedge\varphi_3, \nn \\
\dd\varphi_2&=\,\,\dd\varphi_3\wedge\varphi_1, \nn \\
\dd\varphi_3&=\,\,\dd\varphi_1\wedge\varphi_2, \nn \\
\dd\varphi_4&=0.
\label{mcsl}
\end{align}

\subsection{A digression: Yang-Mills equations}
\label{subsec:YM}
We consider the matrix with 1-form entries
\beq
\mathfrak A\,=\,\alpha\sum_{a=1}^3(i\tau_a\otimes\varphi_a)+\beta(\sigma_0\otimes\varphi_4),
\label{depv}
\eeq
which is a $\tilde{\mathfrak{g}}$-valued $1$-form on $\tilde A$ with $\sigma_{0}=1_2$ the two dimensional identity matrix and the two dimensional matrices $\tau$ given by the Pauli matrices for the $\mathfrak{su}(2)$ case, while  $(\tau_1=\sigma_1, \,\tau_{2,3}=i\sigma_{2,3})$ for $\mathfrak{sl}(2,\R)$ . 

As already mentioned above, since both $\mathfrak{su}(2)$ and $\mathfrak{sl}(2,\R)$ are semisimple, the corresponding Cartan - Killing form $g_{ab}=f_{ak}^{\,\,\,\,s}f_{bs}^{\,\,\,\,k}$ is non degenerate. Therefore the tensor 
(with $\lambda\in\R$)
\beq
\label{gsu2}
\gamma\,=\,\gamma_{\mu\nu}\,\varphi_{\mu}\otimes\varphi_{\nu}\,=\,u_4^2\,g_{ab}\,\varphi_a\otimes\varphi_b+\lambda^2\varphi_4\otimes\varphi_4
\eeq
gives a symmetric non degenerate 2-form on the cotangent space for the calculus on $\tilde A$.  The problem of defining a consistent Hodge duality on the exterior algebra for a non commutative algebra has not a general solution\footnote{When the commutative algebra of functions on two and three dimensional spheres is deformed within the quantum group setting, a Hodge duality on the corresponding exterior algebras has been studied in \cite{post-stellen, alemeron, ale1112, ale11, laza, lareza1}. Within the Drinfeld - Jimbo deformation formalism we refer the reader to \cite{dimit, aschi} and references therein for  a consistent notion of a Hodge duality for a a differential  calculus suitably given by a twist quantization. }.  Given the four dimensional calculi we have developped on $\tilde{A}$, the identities \eqref{eafr} allow to set,  with respect to the volume form $\Omega=(\mathrm{sign}\,g)\lambda u_4^3\,\varphi_1\wedge\varphi_2\wedge\varphi_3\wedge\varphi_4$, the relations
\begin{align}
&*_H \,\varphi_{\mu}\,=\,\tilde{\gamma}_{\mu\sigma}\,i_{D_{\sigma}}\Omega, \nn \\
&*_H\,(\varphi_{\mu_1}\wedge\dots\wedge\varphi_{\mu_k})\,=\,\tilde{\gamma}_{\mu_1\sigma_{a_1}}\dots\tilde{\gamma}_{\mu_k\sigma_{a_k}}\,i_{D_{a_1}}\dotsi_{D_{a_k}}\Omega.
\label{hodgef}
\end{align}
Such relations consistently define a Hodge duality on the frame for the calculus, with $\tilde{\gamma}_{\mu\nu}\gamma_{\nu\rho}=\delta_{\mu\rho}$. The Hodge duality is then extended to the whole exterior algebra by requiring it to be $\tilde{A}$-linear. Given such a Hodge structure, and recalling that $\dd u_4=0$, it is straightforward to prove the following identities, 
\begin{align}
&\mathfrak A\wedge \mathfrak A=-(2\alpha)\dd\mathfrak A \nn \\
&\mathfrak F=\dd\mathfrak A+\mathfrak A\wedge\mathfrak A=(1-2\alpha)\dd\mathfrak A \nn \\
&*_H\mathfrak F=\lambda u_4^{-1}(2\alpha-1)(\mathfrak A\wedge\varphi_4) \nn \\
&\dd(*_H\mathfrak F)\,+\,\mathfrak A\wedge(*_H\mathfrak F)-(*_H\mathfrak F)\wedge \mathfrak A=\lambda u_4^{-1}(2\alpha-1)(2-\frac{1}{2\alpha})(\mathfrak A\wedge\mathfrak A\wedge\varphi_4)
\label{idesu2}
\end{align}
with $*_H$ the Hodge star operator. 
The second line gives the curvature $\mathfrak F$ corresponding to the  vector potential $\mathfrak A$, the fourth line shows that the covariant derivative of $*_H\mathfrak F$ is zero if and only if it is either $\alpha=1/2$ (a case that we discard  since it corresponds to a zero curvature: notice that it comes from the pure gauge term given by the Maurer-Cartan form for $\mathfrak g$) or $\alpha=1/4$. This means that the vector potential $\mathfrak A$ gives a non trivial  solution of the homogeneous (i.e. sourceless) Yang-Mills equations for $\alpha=1/4$. The above calculations show that the existence of this solution, which is valid for any value of the parameters $\beta$ and $\lambda\neq0$ in the metric tensor $g$, only depends on the semisimplicity of the Lie algebra $\mathfrak{su}(2)$. This solution is the analogue on the noncommutative space $\tilde A$ of the meron solution for the Yang-Mills equations which was introduced in \cite{daff}, analysed within the principal and vector bundle formalism in \cite{lame}, and extended to the quantum group setting in \cite{alemeron}. 

We close this digression upon noticing that, for the algebra $\tilde{A}$ corresponding to the $\mathfrak{su}(2)$ case, the Hodge - de Rham Laplacian operator given by 
$$
\Box\,=\,*_H\dd*_H\dd\,=\,\delta_{ab}[u_a,[u_b, ~]_{\theta}]_{\theta}
$$ 
coincides with the well studied Laplacian for the fuzzy sphere, whose spectral dimension is equal to 2. The problem of defining a Laplacian with spectral dimension equal to 3 can be studied within the fuzzy space approach, and an interesting example along this path is provided by \cite{galpre1,galpre2}.

\section{A differential calculus on $A_{\mathfrak g}$ for a non semisimple Lie algebra $\mathfrak g$}
\label{sec:not}

We begin this section upon recalling that, within the classical formalism of differential geometry, for a  Lie algebra $\mathfrak g$ given by $[x_a,x_b]=f_{ab}^{\,\,\,\,c}x_c$ there exists a natural Poisson structure on $\mathfrak g^*$ given by (notice that we are identifying the coordinates on $\mathfrak g$ with those on $\mathfrak g^*$)
$$
\Lambda_{\mathfrak g}\,=\,\frac{1}{2}\,f_{ab}^{\,\,\,\,c}x_c\,\frac{\partial}{\del x_a}\,\wedge\,\frac{\partial}{\partial x_b}.
$$
A vector field $V\,=\,B_{ab}x_a\del_b$ is a derivation for $\mathfrak g$ if and only if it satisfies the relation
$$
L_V\,\Lambda_{\mathfrak g}=0.
$$
For a non semisimple $\mathfrak g$ there exist vector fields $V$ whose action upon the corresponding  universal envelopping algebra cannot be represented in terms of a commutator with an element in $\mathfrak g$. Such operators are naturally called \emph{outer} derivations on $\mathfrak g$. 

In this section we describe how the union of inner and outer derivations for a non semisimple Lie algebra $\mathfrak g$ allows to define 
a tangent space for a differential calculus on $A_{\mathfrak g}$ which turns out to be parallelisable.

\subsection{The case $\mathfrak g=\mathfrak e(2)$}
\label{subsec:e2}
The Lie algebra $\mathfrak g=\mathfrak e(2)$ is given in \eqref{je2}, so we have 
$$
\Lambda_{\mathfrak e(2)}\,=\,
\frac{1}{2}\,(x_3\,\del_1\wedge\del_2\,+\,x_2\del_3\wedge\del_1).
$$
The vector fields 
\begin{align}
&V_1\,=\,x_3\del_2-x_2\del_3, \nn \\ &V_2\,=\,-x_3\del_1, \nn \\ &V_3\,=\,x_2\del_1 \label{inde2}
\end{align}
give the inner derivations for $\mathfrak e(2)$, with $V_k(~)\,=\,[x_k,~]$, while the vector field
\beq
\label{ede2}
V_E\,=\,x_2\del_2+x_3\del_3
\eeq
is an exterior derivation for the universal envelopping algebra corresponding to $\mathfrak e(2)$. From the relations 
\begin{align}
& [V_1,V_E]\,=\,0,  \nn \\ & [V_2,V_E]\,=\,-V_2, \nn \\  &[V_3,V_E]\,=\,-V_3 \label{gte2}
\end{align}
we see that the set $\{V_{\mu}\}_{\mu=1,\dots,4}=\{V_1, V_2, V_3, V_E\}$ gives a basis for a Lie algebra $\tilde{\mathfrak e}(2)$ which is a one dimensional extension of  $\mathfrak e(2)$. We  write the Lie algebra structure for $\tilde{\mathfrak e}(2)$ as 
\beq
\label{tie2}
[V_{\mu},V_{\nu}]\,=\,\tilde{f}_{\mu\nu}^{\,\,\,\,\,\rho}V_{\rho}.
\eeq
We define  the noncommutative algebra $A_{{\mathfrak e}(2)}\subset\mathcal M^{\theta}$ in terms of  the elements $(u_1,u_2,u_3)$ given in \eqref{je2} as the quotient \eqref{Aqu}.  Notice that this algebra coincides with $\tilde A$  defined by 
$$
\tilde A\,=\,\{f\,\in\,\mathcal M^{\theta}\,:\,[u_C,f]_{\theta}=0\},
$$
with $u_C\,=\,\frac{1}{2}(u_2^2+u_3^2)\,=\,\frac{1}{2}(q_1^2+q_2^2)$. 
It is clear that the operators
\beq
\label{inqe2}
D_k(f)\,=\,[u_k,f]_{\theta}
\eeq
define inner derivations on $\tilde{A}=A_{{\mathfrak e}(2)}$. The operator $D_E$, whose action upon  $\tilde{A}$ is defined by setting 
\begin{align}
D_E&:\,u_1\,\mapsto\,0, \nn \\
&:u_2\,\mapsto\,(i\theta)u_2, \nn \\
&:u_3\,\mapsto\,(i\theta)u_3
\label{de2DE} 
\end{align}
on the generators and extended in terms of the Leibniz rule, is an outer derivation on $\tilde A$. This action can also be written as 
\beq
\label{De2i4}
D_E(f)\,=\,[u_E,f]_{\theta} \qquad \mathrm{for}\,\,f\,\in\,\tilde A
\eeq
with 
\beq
\label{duEe2}
\mathcal M^{\theta}\,\ni\,u_E\,=\,-(q_1p_1+q_2p_2).
\eeq 
This means that the action of the outer derivation $D_E$ for $\tilde{A}$ can be represented as a commutator on $\tilde A\subset \mathcal M^{\theta}$ in terms of a quadratic element $u_E\,\in\,S\,\subset\,\mathcal M^{\theta}$. 
Notice that the element $u_E$ is defined up to an arbitrary function of the quadratic Casimir $u_C$, but this does not affect any of the results we shall describe.

The set $\{u_{\mu}\}_{\mu= 1,\dots,4}\,=\,\{u_1,u_2,u_3,u_E\}$ is  a Poisson subalgebra in $(\mathcal M^{\theta}, *)$ which is isomorphic (up to the factor $i\theta$) to the Lie algebra $\tilde{\mathfrak e}(2)$ written in \eqref{tie2}. We denote by $\mathcal D\,\simeq\,\tilde{\mathfrak e}(2)$ the Lie algebra of derivations spanned by $\{D_{\mu}\}_{\mu=1,\dots,4}=\{D_k,D_E\}$  and describe the differential calculus on $\tilde{A}$ based on it. 

Since $\tilde{\mathfrak e}(2)\,\subset\,\mathfrak{isp}(4,\R)$, we consider the elements  $\{\alpha_{\mu}\}_{\mu=1,\dots,4} \,\in\,
(C_{\wedge}(\mathfrak{isp}(4,\R), \mathcal M^{\theta}), \dd)$
whose action is given by  
\beq
\label{fofae2}
\alpha_{\mu}(D_{\nu})\,=\,[u_{\nu},u_{\mu}]_{\theta}
\eeq
and clearly observe that, on $\mathcal M^{\theta}$, it is  $\alpha_{\mu}=\dd u_{\mu}$. The elements in $C^1_{\wedge}(\mathfrak{isp}(4,\R), \mathcal M^{\theta})$ given by 
\begin{align}
&\omega_1\,=\,\frac{1}{2}(u_3\alpha_2-u_2\alpha_3), \nn \\
&\omega_2\,=\,-\frac{1}{2}(u_3\alpha_1+u_2\alpha_E), \nn \\
&\omega_3\,=\,\frac{1}{2}(u_2\alpha_1-u_3\alpha_E), \nn \\
&\omega_E\,=\,\frac{1}{2}(u_2\alpha_2+u_3\alpha_3)
\label{ode2}
\end{align} 
verify (with $\mu=1,\dots,4$)
\beq
\label{duae2om}
\omega_{\mu}(D_{\sigma})\,=\,(i\theta)u_C\delta_{\mu\sigma},
\eeq
so the elements 
\beq
\label{basdue2}
\varphi_{\mu}\,=\,-\frac{i}{\theta}\,u_C^{-1}\omega_{\mu}
\eeq
give a basis for $C^1_{\wedge}(\tilde{\mathfrak e}(2), A_{\mathfrak e(2)})$ (where we have  -- analogously to the localisation considered  in the previous sections --  extended the algebra $A_{\mathfrak e(2)}$ upon adding\footnote{Notice that the algebra $A_{\mathfrak e(2)}$ given by such a localisation can be seen as a noncommutative deformation of the algebra $\mathcal F(\R^3\backslash\{u^2_2+u_3^2=0\})$, i.e. of the functions defined on $\R^3$ without a straight line.} the term $u_C^{-1}$ with $u_C^{-1}u_k=u_{k}u_C^{-1}$ and $u_Cu_C^{-1}=u_C^{-1}u_C=1$) which is dual to the basis $\{D_{\mu}\}_{\mu=1,\dots,4}$. Using the \eqref{wedp} for a wedge product, the definition \eqref{ddef} to introduce a $\dd$ operator and the definition \eqref{contX} for the contraction operator, we have that $(C_{\wedge}(\tilde{\mathfrak e}(2), A_{\tilde{\mathfrak e}(2)}), \dd)$ provides a four dimensional differential calculus on $A_{\tilde{\mathfrak e}(2)}$. As for the calculi described in the previous section (see \eqref{eafr}), the identities $f*\varphi_{\mu}\,=\,\varphi_{\mu}*f$ and 
$\varphi_{\mu}\wedge\varphi_{\sigma}=-\varphi_{\sigma}\wedge\varphi_{\mu}$ are valid. 
It is therefore  immediate to prove that the action of the exterior derivative 
$\dd\,:\,A_{\tilde{\mathfrak e}(2)}\,\to\,C^1_{\wedge}(\tilde{\mathfrak e}(2), A_{\tilde{\mathfrak e}(2)})$ 
can be written as
\beq
\dd f\,=\,(D_{\mu}f)\varphi_{\mu}\,=\,([u_{\mu},f]_{\theta})\varphi_{\mu}
\label{free2}
\eeq
in terms of the commutator structure on $\mathcal M^{\theta}$. 
The Maurer - Cartan relation for such a calculus comes from the Lie algebra structure of $\tilde{\mathfrak e}(2)$ as in \eqref{tie2}, 
$$
\dd\varphi_{\rho}\,=\,-\frac{1}{2}\,i\theta\tilde f_{\mu\nu}^{\,\,\,\,\,\rho}\varphi_{\mu}\wedge\varphi_{\nu}
$$
and reads explicitly as follows
\begin{align}
&\dd\varphi_1\,=\,0, \\
&\dd\varphi_2\,=\,i\theta(\varphi_3\wedge\varphi_1+\varphi_2\wedge\varphi_4), \nn \\
&\dd\varphi_3\,=\,i\theta(\varphi_2\wedge\varphi_1+\varphi_3\wedge\varphi_4), \nn \\
&\dd\varphi_4\,=\,0. 
\label{mcee2}
\end{align}
One can also immediately see that the presence of a 1-form which dualises the outer derivation for $A_{\mathfrak e(2)}$ allows to prove that centre $Z(A_{{\mathfrak e}(2)})$ of the algebra is not in the kernel of the exterior derivative operator $\dd$ introduced in this section. This result comes immediately upon computing that 
\beq
\label{dcae2}
\dd u_C\,=\,2(i\theta)u_C\,\varphi_4.
\eeq

\subsection{The case $\mathfrak g=\mathfrak{iso}(1,1)$}
\label{subsec:iso11}
The Lie algebra $\mathfrak{iso}(1,1)$ comes as a real form of the complex Lie algebra $~$. It follows that the elements
\begin{align}
f_1=iu_1, \nn \\
f_2=u_2, \nn \\
f_3=iu_3, \label{compmap}
\end{align}
(where $u_{k=1,\dots,3}\in\,\mathcal S$ are the generators of the Lie algebra ${\mathfrak e}(2)$, see \eqref{je2}), give a realisation of 
$A_{\mathfrak{iso}(1,1)}$. Under the identification \eqref{compmap} one sees that the analysis performed in the previous subsection for the $A_{\mathfrak e(2)}$ algebra can be -- \emph{mutatis mutandis} -- carried through for the $A_{\mathfrak{iso}(1)}$ case. Along this path we can then equip $A_{\mathfrak{iso}(1,1)}$ with a derivation based four dimensional differential calculus. 
  
\subsection{The case $\mathfrak g=\mathfrak h(1)$}
\label{subsec:h1}
The Lie algebra $\mathfrak{g}=\mathfrak h(1)$ corresponding to the Heisenberg-Weyl group is given in \eqref{jh1}. The corresponding Poisson structure is given by 
\beq
\label{h1poi}
\Lambda\,=\,x_3\,\frac{\del}{\del x_1}\,\wedge\,\frac{\del}{\del x_2}.
\eeq
For such a Lie algebra there are two inner derivations, whose action is given by the vector fields 
\begin{align}
&V_1\,=\,x_3\del_2, \nn \\
&V_2\,=\,-x_3\del_1, 
\label{indh1}
\end{align}
and four exterior derivations, whose action is given by the  vector fields 
\begin{align}
&E_1\,=\,x_1\del_1+x_3\del_3, \\
&E_2\,=\,x_2\del_2+x_3\del_3, \\
&E_R\,=\,x_1\del_2-x_2\del_1, \\
&E_H\,=\,x_1\del_2+x_2\del_1.
\label{estdh1}
\end{align}
The commutator structure closed by these elements reads interesting Lie subalgebras. We have
\begin{align}
&[V_1,V_2]=0, \nn \\
&[E_1,E_2]=0,  \nn \\
&[V_a,E_b]=-\delta_{ab}V_a 
\label{suh1}
\end{align}
for $a,b=1,2$, with 
\begin{align}
[E_R,V_1]=-V_2, &\qquad\qquad [E_R,V_2]=V_1 \nn \\
[E_H,V_1]=V_2, & \qquad\qquad [E_H,V_2]=V_1
\label{erhV}
\end{align}
and 
\begin{align}
&[E_R,E_H]=2(E_1-E_2), \nn \\
&[E_R,E_1]=-E_H,\qquad\qquad[E_R,E_2]=E_H, \nn \\
&[E_H,E_1]=-E_R,\qquad\qquad[E_H,E_2]=E_R.
\label{erh12}
\end{align}
These relations show that $\{V_a,E_b\}_{a,b=1,2}$ span a Lie algebra $\tilde{\mathfrak g}$ that extends the abelian Lie algebra spanned by $V_1,V_2$. We see that also the outer derivations  $\{E_1,E_2, E_R,E_H\}$ alone span a Lie algebra. 

Within the noncommutative setting, we consider $A_{\mathfrak h(1)}\subset \mathcal M^{\theta}$ to be the algebra generated (and suitably completed) as in \eqref{Aqu} by the elements $u_1,u_2,u_3$ defined in \eqref{jh1}. This algebra coincides with the algebra $\tilde{A}$ defined by
$$
\tilde{A}\,=\,\{f\,\in\,\mathcal M^{\theta}\,:\,[u_3,f]_{\theta}\,=\,0\}.
$$
The operators 
\beq
\label{inqh1}
D_k(f)\,=\,[u_k,f]_{\theta}
\eeq
define inner derivations on $A_{\mathfrak h(1)}$. Analogously to what we described for the $\mathfrak g=\mathfrak e(2)$ case, the operators $D_{E_a}$ with $a=1,2$ whose action upon $A_{\mathfrak h(1)}$ is defined by setting
\beq
\begin{array}{lcl} D_{E_1}(u_1)\,=\,(i\theta)u_1, &\qquad\qquad & D_{E_2}(u_1)\,=\,0,  \\
D_{E_1}(u_2)\,=\,0, &\qquad\qquad & D_{E_2}(u_2)\,=\,(i\theta)u_2,  \\
D_{E_1}(u_3)\,=\,(i\theta)u_3, &\qquad\qquad &  D_{E_2}(u_3)\,=\,(i\theta)u_3 \end{array}
\label{deeh1}
\eeq
on the generators and extended in terms of the Leibniz rule, are outer derivations for $A_{\mathfrak h(1)}$. This action can be written as 
\beq
\label{dech1}
D_{E_a}(f)\,=\,[u_{E_a},f]_{\theta}
\eeq
with 
\begin{align}
&u_{E_1}=-(p_1q_1+p_2q_2), \nn \\ &u_{E_2}=-p_2q_2.
\label{defEh1}
\end{align}
The analogy with the analysis performed in the $\mathfrak e(2)$ case is evident. 
The action of the outer derivation $D_{E_a}$ for $A_{\mathfrak h(1)}$ can be represented as a commutator on $ A_{\mathfrak h(1)}\subset \mathcal M^{\theta}$ in terms of the\footnote{Notice further that the elements $u_{E_a}$ are defined up to an arbitrary function of the element $u_3$, but this does not affect any of the results given in this section.} quadratic element $u_{E_a}\,\in\,S\,\subset\,\mathcal M^{\theta}$.

It is nonetheless easy to prove that there are no elements $u_R,u_H\,\in\,\mathcal M^{\theta}$ such that the   action of the derivations $D_R,D_H$ on $A_{\mathfrak h(1)}$ associated to the vector fields $E_R, E_H$, given by 
\beq
\begin{array}{lcl}
D_{R}(u_1)\,=\,-(i\theta)u_2, &\qquad\qquad &D_{H}(u_1)\,=\,(i\theta)u_2, \\
D_{R}(u_2)\,=\,(i\theta)u_1, &\qquad\qquad &D_{H}(u_2)\,=\,(i\theta)u_1,  \\
D_{R}(u_3)\,=\,0, &\qquad\qquad &D_{H}(u_3)\,=\,0
\end{array}
\label{derhh1}
\eeq
on the generators, can be written as a commutator with $u_R,u_H$ as in \eqref{dech1}.

In order to introduce a derivations based calculus on $A_{\mathfrak h(1)}$ we consider the elements $\{u_{\sigma}\}_{\sigma=1,\dots,4}=\{u_1,u_2,u_3,u_4=-\mu p_1q_1-\nu p_2q_2\}$, with $\mu, \nu\,\in\,\R$. These elements span, with respect to the non commutative product,  a family (depending on the parameters $\mu, \nu$) of Lie subalgebras in $\mathcal M^{\theta}$ which are isomorphic (up to the factor $(i\theta)$)  to the four dimensional  Lie algebras $\tilde{\mathfrak h}(1)$ which extend $\mathfrak h(1)$ as they are  spanned by the elements  
$$
\tilde{\mathfrak h}(1)\,=\,\{\mathfrak h(1),\,\mu E_1+(\nu-\mu)E_2\}.
$$
The relevant Lie  algebra structure is given by 
\begin{align}
&[u_1,u_2]_{\theta}=(i\theta)u_3, \nn \\
&[u_1,u_4]_{\theta}=-(i\theta)\mu u_1, \nn \\
&[u_2,u_4]_{\theta}=(i\theta)(\mu-\nu)u_2, \nn \\
&[u_3,u_4]_{\theta}\,=\,-(i\theta)\nu u_3.
\label{h1tcs}
\end{align}
In analogy to what we wrote in the previous subsection, we describe the differential calculus on $A_{\mathfrak h(1)}$ based on the Lie algebra $\tilde{\mathfrak h}(1)$  derivations $D_{\sigma}$ given by
\beq
D_{\sigma}(f)\,=\,[u_{\sigma},f]_{\theta}.
\label{famDmn}
\eeq
Since $\tilde{\mathfrak h}(1)\,\subset\,\mathfrak{isp}(4,\R)$, we consider the elements 
$\{\alpha_{\rho}\}_{\rho=1,\dots,4} \,\in\,
(C_{\wedge}(\mathfrak{isp}(4,\R), \mathcal M^{\theta}), \dd)$
whose action is given by  
\beq
\label{foht1}
\alpha_{\rho}(D_{\sigma})\,=\,[u_{\sigma},u_{\rho}]_{\theta}
\eeq
and observe that $\alpha_{\rho}\,=\,\dd u_{\rho}$. It is a long but straightforward calculation to see that, if $\nu\neq0$, then the elements 
\begin{align}
&\omega_1\,=\,\frac{1}{2}\,(u_3\alpha_2\,+\,(\frac{\mu}{\nu}\,-1)u_2\alpha_3), \nn \\
&\omega_2\,=\,\frac{1}{2}\,(-u_3\alpha_1\,+\,\frac{\mu}{\nu}\,u_1\alpha_3), \nn \\
&\omega_3\,=\,\frac{1}{2}\,((1\,-\,\frac{\mu}{\nu})u_2\alpha_1\,-\,\frac{\mu}{\nu}\,u_1\alpha_2\,+\,(i\theta)(\frac{\mu}{\nu}\,-\,\frac{\mu^2}{\nu^2})\alpha_3\,-\,\frac{1}{\nu}\,u_3\alpha_4), \nn \\
&\omega_4\,=\,\frac{1}{2\nu}\,u_3\alpha_3
\label{framht1}
\end{align}
verify 
\beq
\label{fradeht1}
\omega_{\rho}(D_{\sigma})\,=\,(i\theta)u_3\delta_{\rho\sigma}.
\eeq
If we denote still by $A_{\mathfrak h(1)}$ the algebra extended  by the localization defined by adding the 
generator $u_3^{-1}$ corresponding to the Casimir function,  with $u_3^{-1}u_3=u_3u_3^{-1}=1$ and  $u_3u_{j}=u_{j}u_3$ for $j=1,2$, then\footnote{Notice that such an extended algebra $A_{\mathfrak h(1)}$  gives  a noncommutative deformation of the classical algebra $\mathcal F(\R^3\backslash\{x_3=0\})$, or more generally of the functions defined on $\R^3$ without a plane.} the elements 
\beq
\label{baduht1}
\varphi_{\rho}\,=\,-\frac{i}{\theta}\,u_3^{-1}\omega_{\rho}
\eeq
give a basis for $C^1_{\wedge}(\tilde{\mathfrak h}(1), A_{\mathfrak h(1)})$. For the corresponding differential calculus  one has that 
$\dd\,:\,A_{\tilde{\mathfrak e}(2)}\,\to\,C^1_{\wedge}(\tilde{\mathfrak e}(2), A_{\tilde{\mathfrak e}(2)})$ 
can be written as
\beq
\dd f\,=\,(D_{\rho}f)\varphi_{\rho}\,=\,([u_{\rho},f]_{\theta})\varphi_{\rho}.
\label{free2}
\eeq
Also for this calculus, as in \eqref{eafr}, it is $f*\varphi_{\mu}\,=\,\varphi_{\mu}*f$ and 
$\varphi_{\mu}\wedge\varphi_{\sigma}=-\varphi_{\sigma}\wedge\varphi_{\mu}$.
The Maurer - Cartan relation  comes from the Lie algebra structure \eqref{h1tcs}  of $\tilde{\mathfrak h}(1)$; as expected, the differential calculus for the Lie algebras $\tilde{\mathfrak h}(1)$  depends on the ratio $\mu/\nu$. It is immediate to see that the identity $\alpha_3=\dd u_3$  given by \eqref{foht1} can be cast as 
\beq
\label{dacca}
\dd u_3\,=\,(i\theta)\nu u_3\,\varphi_4,
\eeq
thus proving that the centre $Z(A_{\mathfrak h(1)})$ of the algebra is not in the kernel of the exterior derivative $\dd$ defined in this section.

\section{Conclusion}\label{concl}
At the end of the paper it is clear that the three dimensional (classical) spaces we mention in the title are given as the foliations of the codimension one regular orbits for the action of the Lie algebra $\mathfrak g$ upon $\mathfrak g^*\simeq\R^3$. For the given noncommutative deformation of  such spaces we have described a derivation based differential calculus and exhibited a frame for it. The natural step forward along the path we started is to analyse  the spinor structures on such exterior algebras. We aim to develop such analysis in a future work.  

\section*{Acknowledgments}
P.V. acknowledges support by COST (European Cooperation in Science
and Technology) in the framework of COST Action MP1405 QSPACE and by  INFN under the project  Iniziativa Specifica GeoSymQFT. G. M.  would like to thank the support provided by the Santander/UC3M
Excellence Chair Programme.


\begin{thebibliography}{99}

\bibitem{aschi} P. Aschieri, A. Borowiec, A. Pacho\l{}, {\it Observables and dispersion relations in $\kappa$-Minkowski spacetime}, JHEP 10(2017)152;


\bibitem{twist} 
  P.~Aschieri, M.~Dimitrijevic, P.~Kulish, F.~Lizzi and J.~Wess,
{\it Noncommutative spacetimes: Symmetries in noncommutative geometry and field theory,}
  Lect.\ Notes Phys.\  774  (2009) 1;
  
\bibitem{aschieri} P.~Aschieri, F.~Lizzi and P.~Vitale,
 {\it Twisting all the way: From Classical Mechanics to Quantum Fields,}
  Phys.\ Rev.\ D  77 (2008)  025037; 
 
\bibitem{oxford} T. Brzezi\'nski, S. Majid, {\it Quantum group gauge theory on quantum spaces}, 157 (1993) 591-638; 
 
  
\bibitem{cmw11} E. Cagnache, T. Masson, J.C. Wallet, {\it noncommutative Yang-Mills-Higgs actions from derivation based differential calculus}, J. Geom. Phys. 5 (2011) 39;

\bibitem{castellani} L. Castellani, {\it Gauge theories of quantum groups}, Phys. Lett. B 292 (1992) 93-98; 


\bibitem{connes} A.~H.~Chamseddine and A.~Connes,
{\it The Spectral action principle,}
  Commun.\ Math.\ Phys.\  186 (1997) 731:



\bibitem{damv} A. D'Avanzo, G. Marmo, A. Valentino, {\it Reduction and unfolding for quantum systems: the hydrogen atom}, Int. J. Geom. Meth. Mod. Phys. 2, 6 (2005) 1043;

\bibitem{daff} V. de Alfaro, S. Fubini, G. Furlan, {\it New classical solutions of the Yang-Mills field equations}, Phys. Lett. B 65 (1976) 163;

  
\bibitem{degmw} A. de Goursac, T. Masson and J.-C. Wallet,
{\it Noncommutative $\epsilon$-graded connections},
J. Noncom. Geom. 6 (2012) 343;


\bibitem{post-stellen} F. Di Cosmo, G. Marmo, J.M. Perez-Pardo, A. Zampini, {\it A Hodge - de Rham Dirac operator on the quantum ${\rm SU}(2)$}, Int. J. of Geom. Meth. Mod. Phys. 15 (2018) 1850030; 

\bibitem{dimit} M. Dimitrijevic, L. Jonke, A. Pacho\l{}, {\it Gauge theory on twisted $\kappa$-Minkowski: old problems and possible solutions}, SIGMA 10 (2014) 063; 

\bibitem{dv01} M. Dubois-Violette, {Lectures on graded differential algebras and noncommutative geometry},  in Maeda et al. (eds),  {\it noncommutative differential geometry and its applications to physics}, Kluwer 2001;


\bibitem{galpre1} V. G\'alikov\'a, S. Kov\'a\v{c}ik, P. Pre\v{s}najder, {\it Laplace-Runge-Lenz vector in quantum mechanics in 
noncommutative spaces}, J. Math. Phys 54 (2013) 122106;


\bibitem{galpre2} V. G\'alikov\'a,  P. Pre\v{s}najder, {\it Coulomb problem in noncommutative quantum mechanics}, J. Math. Phys. 54 (2013) 052102;


\bibitem{marse} V. Gayral, J.M. Gracia-Bond\'ia, B. Iochum, T. Sch\"ucker, J.C. Varilly, {\it Moyal planes are spectral triples}, Comm. Math. Phys. 246 (2004) 569;

\bibitem{gjw15} A. G\'er\'e, T. Juri\'c, J.C. Wallet, {\it Noncommutative gauge theories on $\R^3_{\lambda}$: perturbatively finite models}, JHEP 12 (2015) 045;

\bibitem{gvw14} A. G\'er\'e, P. Vitale, J.C. Wallet, {\it Quantum gauge theories on noncommutative three-dimensional spaces}, Phys. Rev. D 90 (2014) 045019;

\bibitem{poipo} J. Grabowski, G. Marmo, A. M. Perelomov, {\it Poisson structures: towards a classification}, Mod. Phys. Lett. A8 (1993) 1719;

\bibitem{sele} J. M. Gracia-Bond\'ia, F. Lizzi, G. Marmo, P. Vitale, {\it Infinitely many star products to play with}, JHEP 04 (2002) 026;

\bibitem{alg1} J. M. Gracia-Bond\'ia, J. C. Varilly, {\it Algebras of distributions suitable for phase space quantum mechanics I}, J. Math. Phys. 29 (1988) 869-879;


\bibitem{gutt83} S. Gutt, {\it An explicit $*$-product on the cotangent bundle of a Lie group}, Lett. Math. Phys. 7 (1983) 249;


\bibitem{shei} A. B. Hammou, M. Lagraa, M. M. Sheikh-Jabbari, {\it Coherent state induced star-product on $\R^3_{\lambda}$ and the fuzzy sphere},  Phys. Rev. D 66 (2002) 025025;

\bibitem{istvan} I. Heckenberger, {\it Spin geometry on quantum groups via covariant differential calculi}, Adv. in Math. 175 (2003) 197-242;


\bibitem{HTS80}
J.\ Harnad, J.\ Tafel and S.\ Shnider, \textit{Canonical connections on Riemannian symmetric spaces and solutions to the Einstein-Yang-Mills equations}, J.\ Math.\ Phys.\ 21 (1980)  2236;

\bibitem{kuvi15} V.G. Kupriyanov, P. Vitale, {\it Noncommutative $\R^d$ via closed star product}, JHEP 08 (2015) 024;

\bibitem{lame} G. Landi, {\it Spinor and gauge connections over oriented spheres}, Proc. 7th Italian Conference on General Relativity and Gravitational Physics, World Scientific (1987) 287;

\bibitem{lama90} G. Landi, G. Marmo, {\it Algebraic differential calculus for gauge theories}, Nucl. Phys. Proc. Suppl. 18A (1990) 171;

\bibitem{lareza1} G.Landi, C.Reina, A.Zampini, {\it Gauged Laplacians on a quantum Hopf
bundle}, Comm. Math. Phys. 287 (2009) 179-209; 

\bibitem{laza} G.Landi, A.Zampini, {\it Calculi, Hodge operators and Laplacians on a quantum Hopf fibration}, 
Rev. Math. Phys. 23 (2011) 575-613; 

\bibitem{matrixLV} F.~Lizzi and P.~Vitale,
{\it Matrix Bases for Star Products: a Review},
  SIGMA  10 (2014) 086 
  
\bibitem{Lizzireview} F. Lizzi, \textit{Noncommutative Geometry and Particle Physics}, Proceedings of the Corfu Summer Institute 2017, arXiv:180500411 [hep-th];

\bibitem{fs} J. Madore, {\it The fuzzy sphere}, Class. Quant. Grav. 9 (1992) 69;



\bibitem{jsnap} V. I. Man'ko, G. Marmo, P. Vitale, F. Zaccaria, {\it A generalization of the Jordan - Schwinger map: classical version and its $q$-deformation}, Int. J. Mod. Phys. A9 (1994) 5541;

\bibitem{dica} G. Marmo, P. Vitale, A. Zampini, {\it Noncommutative differential calculus for Moyal subalgebras}, J. of Geom.  Phys. 56 (2006) 611;

\bibitem{masson08} T. Masson, {\it Examples of derivation based differential calculi related to noncommutative gauge theories}, Int. J. Geom. Meth. Mod. Phys. 05 (2008) 1315; 


\bibitem{meli1} S. Meljanac, Z. \v{S}koda, M. Stojic, {\it Lie algebra type non commutative phase spaces are Hopf algebroids}, Lett. Math. Phys. 107 (2017) 475-503; 



\bibitem{meli2} S. Meljanac, S. Kre\v{s}ic-Juric, T. Martinic, {\it Realization of bicovariant differential calculus on the Lie algebra type non commutative spaces}, J. Math. Phys. 58 (2017) 071701;










\bibitem{NT78} J. Nowakowski, A. Trautman, \textit{Natural connections on Stiefel bundles are sourceless gauge fields}, J. Math. Phys.  19 (1978) 1100;

\bibitem{rovi12} L. Rosa, P. Vitale, {\it On the $*$-product quantization and the Duflo map in three dimensions}, Mod. Phys. Lett. A27 (2012) 1250207;

\bibitem{segal67} I.E. Segal, {\it Quantized differential forms}, Topology 8 (1967) 147;

\bibitem{segal70} I.E. Segal, {\it Quantization of the de Rham complex}, Proc. Symp. Pure Math. 16 (1970) 205;

\bibitem{vitale14} P. Vitale, {\it Noncommutative field theory on $\R^3_{\lambda}$}, Fort. Phys. 62 (2014) 825;

\bibitem{alg2} J. C. Varilly, J. M. Gracia-Bond\'ia, {\it Algebras of distributions suitable for phase space quantum mechanics II -- Topologies on the Moyal algebra}, J. Math. Phys. 29 (1988) 880-887; 

\bibitem{viwa13} P. Vitale, J.C. Wallet, {\it Noncommutative field theories on $\R^3_{\lambda}$: towards UV/IR mixing freedom}, JHEP 04 (2013) 115;


\bibitem{wallet09}  J.C. Wallet, {\it Derivations of the Moyal algebra  and noncommutative gauge theories}, SIGMA 5 (2009) 013;

\bibitem{woro1} S.L. Woronowicz, {\it Twisted $\rm SU(2)$ group. An example of a non commutative differential calculus}, 23 (1987) 117-181;

\bibitem{woro2} S.L. Woronowicz, {\it Differential calculus on compact matrix pseudogroups (quantum groups)}, Comm. Math. Phys. 122 (1989) 125-170;


\bibitem{ale1112} A.Zampini, {\it Hodge duality operators on left covariant exterior algebras over two and three dimensional  quantum  spheres}, Rev. Math. Phys. 25 (2013) 9-38; 

\bibitem{ale11} A.Zampini, {\it (A class of) Hodge duality operators over the quantum ${\rm SU(2)}$}, 
J. Geom. Phys. 62 (2012) 1732-1746;

\bibitem{alemeron} A. Zampini, {\it Warped products and Yang-Mills equations on noncommutative spaces}, Lett. Math. Phys. 105 (2015) 221;




\end{thebibliography}
\end{document}